\documentclass[11pt]{amsart} 
\usepackage[T1]{fontenc} \usepackage[utf8]{inputenc} 
\usepackage{lmodern} 
\usepackage{amsmath,amssymb,amsfonts,amsthm} 
\usepackage{dirtytalk}
\usepackage{mathtools}
\usepackage{enumitem} 
\usepackage{graphicx}
\usepackage{xcolor} 
\usepackage{hyperref}
\usepackage[nameinlink]{cleveref} 
\newtheorem{theorem}{Theorem}[section] 
\newtheorem{proposition}[theorem]{Proposition} 
\newtheorem{lemma}[theorem]{Lemma} 
\newtheorem{corollary}[theorem]{Corollary} 
\theoremstyle{definition}

\newtheorem{remark}[theorem]{Remark} 
\newcommand{\R}{\mathbb{R}} 
\newcommand{\C}{\mathbb{C}} 
\newcommand{\N}{\mathbb{N}} 
\newcommand{\Z}{\mathbb{Z}} 
\newcommand{\Lp}[1]{L^{#1}} 
 
\DeclareMathOperator{\supp}{supp} 
\DeclareMathOperator{\Op}{Op}

\title{$L^p$-boundedness of magnetic pseudodifferential operators}

\author{Morten Nielsen} 
\email{mnielsen@math.aau.dk}
\author{Gabriele Peluso}
\email{gpe@math.aau.dk}
\author{Mikkel Hviid Thorn}
\email{mikkelht@math.aau.dk}
\address{Department of Mathematical Sciences, Aalborg University} 
 
\begin{document} 

\begin{abstract}

We show that $L^p$-boundedness for magnetic pseudodifferential operators is equivalent to $L^p$-boundedness of non-magnetic pseudodifferential operators for $p\in[1,\infty]$, regular magnetic fields and classical H\"ormander symbols $S_{\rho,\delta}^m$ with $0<\rho\leq 1$ and $0\leq\delta<1$.

As a corollary we get $L^p$-boundedness for $1<p<\infty$ of magnetic pseudodifferential operators when $\rho>0,0\leq\delta\leq\rho\leq1,\delta<1$ and $m\leq-d(1-\rho)|1/2-1/p|$. Moreover, we give a bound on the operator norm in terms of the magnetic field and symbol. 

\end{abstract} 
\maketitle 
\section{Introduction} 

Continuity of pseudodifferential operators has been of interest since the inception of the pseudodifferential calculus with e.g. the much celebrated Calder\'on--Vaillancourt Theorem \cite{Calderon1972} proving sharp results on the $L^2$-boundedness of pseudodifferential operators. Fefferman \cite{Fefferman1973} gave sharp results of $L^p$-boundedness of pseudodifferential operators with non-multiplier H\"ormander symbols, which was improved by P\"aiv\"arinta and Somersalo \cite{Somersalo1988}. See \cite[pp. 322--323]{Stein2016} for additional history and sources.

For the magnetic pseudodifferential calculus $L^2$-boundedness (and boundedness in magnetic $L^2$-Sobolev spaces), i.e. a magnetic Calder\'on--Vaillancourt Theorem, was first shown in \cite{Purice2007} for the magnetic Weyl quantization and a variant for magnetic Fourier integral operators was shown in \cite{IftimiePurice2011}. Later matrix decompositions \cite{CorneanGarde2019} and frame techniques \cite{Cornean2024,Thorn2026} were applied, where \cite{CorneanGarde2019,Thorn2026} worked with symbol classes beyond the classical H\"ormander classes and \cite{Thorn2026} with asymmetrical magnetic quantizations. Quite recently a preprint \cite{FreyWeng2026} constructed magnetic modulation spaces and proved boundedness results on these.


\subsection{Main results}

Our main theorem is on the equivalence of $L^p$-boundedness of magnetic pseudodifferential operators and of non-magnetic pseudodifferential operators. All definitions and notation used is properly introduced in Section \ref{sec:pre}, in particular the definition of a regular magnetic field is given in Equation \eqref{eq:magdef}. Note that \cite[Lemma 1.2]{Purice2007} and \cite[Proposition 21]{Purice2004} imply that magnetic pseudodifferential operators extend uniquely to continuous maps on tempered distributions. So by $L^p$-boundedness, $p\in[1,\infty]$, we mean that the operators are invariant on $L^p(\R^d)\subseteq\mathcal{S}'(\R^d)$ and bounded on $L^p(\R^d)$.

\begin{theorem}\label{thm:mainresult}
    Let $B$ be a regular magnetic field with potential $A$ with components in $ C_\mathrm{pol}^\infty(\R^d)$, i.e. $dA=B$. Let also $p\in[1,\infty]$, $0<\rho\leq1$, $0\leq\delta<1$ and $m\in\R$.
    
    For $a\in S_{\rho,\delta}^m(\R^{2d})$ the operator $\operatorname{Op}_\mathrm{w}^A(a)$ is bounded on $L^p(\R^d)$ if and only if $\operatorname{Op}_\mathrm{w} (a)$ is bounded on $L^p(\R^d)$.
\end{theorem}

The proof is given in Section \ref{sec:proof}. The strategy behind the proof consists of comparing the distributional kernels $K^A(a)$ and $K^0(a)$ of $\operatorname{Op}^A_\mathrm{w}(a)$ and $\operatorname{Op}_\mathrm{w}(a)$ respectively. They are related by $K^A(a)(x,y)=\Lambda^A(x,y)K^0(a)(x,y)$, i.e. $K^A(a)$ is the \say{magnetic} twisted integral kernel of $K^0(a)$ \cite{Purice2018}. To study it, we split the kernels into two components: one with \say{support} in a strip around the diagonal $|x-y|<R$ and the other having support off-diagonally. Since $\Lambda^A$ is uni-modular, the off-diagonal parts can be treated with a simple application of the Schur test, but this requires the hypothesis $\rho>0$.
The diagonal components require additional care and we further split them into countable components with \say{supports} in cubes near the diagonal. The regularity of the magnetic field allows us to expand the phase $\Lambda^A$ in any cube as a Fourier series with Fourier coefficients that are uniformly dominated by an absolutely convergent series which is independent on the chosen cube. This implies that the operators defined by the integral kernel with cubic supports are uniformly bounded and a vector-valued Schur test concludes the proof.

\begin{remark}
    Note the Theorem also establishes the equivalence of $L^p$-boundedness between magnetic pseudodifferential operators with different magnetic fields as long as they are regular magnetic fields.
    
    Also the arbitrariness of the associated magnetic potential, as long as $A\in C_\mathrm{pol}^\infty(\R^d;\R^d)$, stems from the fact that the magnetic pseudodifferential calculus is gauge-covariant, that is if $dA=d\tilde{A}=B$, then $\tilde{A}=A+\nabla\varphi$ for some $\varphi\in C_\mathrm{pol}^\infty(\R^d)$ and $e^{i\varphi(X)}\mathrm{Op}_\mathrm{w}^A(a)e^{-i\varphi(X)}=\mathrm{Op}_\mathrm{w}^{\tilde{A}}(a)$ in $\mathcal{L}(\mathcal{S}(\R^d),\mathcal{S}'(\R^d))$ with $a\in\mathcal{S}'(\R^{2d})$, see \cite{Purice2004,Purice2007}, where $e^{i\varphi(X)},e^{-i\varphi(X)}$ are $L^p$-isometries. This allows us to work with a definite choice of magnetic potential, and we choose the one in the transversal gauge (to be introduced later).
\end{remark}

We also have a secondary result, combining results on $L^p$-boundedness of non-magnetic pseudodifferential operators with the above result:

\begin{corollary}\label{cor:mainresult}
    Let $B$ be a regular magnetic field with potential $A$ with components in $ C_\mathrm{pol}^\infty(\R^d)$. Let also $p\in(1,\infty)$ and $\rho>0,0\leq\delta\leq\rho\leq1,\delta<1$.
    
    If $m\leq -d(1-\rho)\Big|\frac{1}{2}-\frac{1}{p}\Big|$ and $a\in S_{\rho,\delta}^m(\R^{2d})$, then $\Op_\mathrm{w}^A(a)$ is bounded on $L^p(\R^d)$. Moreover, there exists $C>0$ and $N,M,K\in\N$ independent of $a$ and $B$ such that
    \[\Vert\Op_\mathrm{w}^A(a)\Vert_{\mathcal{L}(L^p)}\leq C(1+\Vert B\Vert_M)^{K}\Vert a\Vert_{S_{\rho,\delta}^m,N}.\]
\end{corollary}

The proof is given in Section \ref{sec:appl}. Note the bound on the operator norm depends on the magnetic field $B$, not the potential $A$, and it also shows that the map $S_{\rho,\delta}^m(\R^{2d})\ni a\mapsto\operatorname{Op}_w^A(a)\in\mathcal{L}(L^p(\R^d))$ is continuous.

\begin{remark}
    Note the results in Corollary \ref{cor:mainresult} are sharp in the sense of \cite{Fefferman1973} since Theorem \ref{thm:mainresult} also, by contraposition, shows equivalence of unboundedness: For every $m>-d(1-\rho)|1/2-1/p|$ there exists $a\in S^m_{\rho, \delta}(\R^{2d})$ such that $\operatorname{Op}^A_\mathrm{w}(a)$ is unbounded on $L^p(\R^d)$. We also remark that Fefferman's counterexample is in fact the  $x$-independent multiplier $e^{i|\xi|^{1-\rho}}\langle\xi\rangle^m$.
\end{remark}

As mentioned, we introduce the different notation and objects in Section \ref{sec:pre}. Section \ref{sec:aux} gives some auxiliary results, necessary for the proof of Theorem \ref{thm:mainresult} in Section \ref{sec:proof}. Next, in Section \ref{sec:appl} we prove Corollary \ref{cor:mainresult}, and lastly, in Section \ref{sec:open} we state some open questions and possible extensions.

\section{Preliminaries}\label{sec:pre}
Let $d\in\N$. We denote by $\mathcal{L}(L^{p_1}(\R^d),L^{p_2}(\R^d))$ the space of bounded operators from $L^{p_1}(\R^d)$ to $L^{p_2}(\R^d)$ and let $\mathcal{L}(L^p(\R^d))\coloneqq\mathcal{L}(L^p(\R^d),L^p(\R^d))$.

\subsection{Pseudodifferential operators}We briefly recall the definition of the classical H\"ormander symbol classes and the associated pseudodifferential operators. 
Set \[ \langle \xi\rangle \coloneqq \bigl(1+|\xi|^2\bigr)^{1/2}, \qquad \xi\in\mathbb R^d. \] 
For multi-indices $j,k\in\mathbb N_0^d$ we use the notation \[ \partial_x^j \coloneqq \partial_{x_1}^{j_1}\cdots \partial_{x_d}^{j_d}, \qquad \partial_\xi^k \coloneqq \partial_{\xi_1}^{k_1}\cdots \partial_{\xi_d}^{k_d}. \] 
Let $m\in\mathbb R$ and $0\leq \delta,\rho\leq 1$. 
The H\"ormander symbol class $S_{\rho,\delta}^m(\mathbb R^d\times\mathbb R^d)$ consists of all functions \[ a\in C^\infty(\mathbb R^d\times\mathbb R^d) \] such that, for every pair of multi-indices $j,k\in\mathbb N_0^d$, there exists a constant $C_{j,k}>0$ for which 
\begin{equation*} 
\bigl| \partial_x^j\partial_\xi^k a(x,\xi) \bigr| \leq C_{j,k} \langle\xi\rangle^{m-\rho|k|+\delta|j|} \quad \forall  (x,\xi)\in\mathbb R^d\times\mathbb R^d,\label{eq:Hormander-symbol-estimate} 
\end{equation*}
and we endow $S_{\rho,\delta}^m(\R^{d}\times\R^d)$ with the topology induced by the semi-norms
\[\Vert a\Vert_{S_{\rho,\delta}^m,N}\coloneqq\sum_{j,k\in\N_0^{d},|j|+|k|\leq N}\sup_{(x,\xi)\in\R^{d}\times\R^d}\langle\xi\rangle^{-m+\rho|k|-\delta|j|}\bigl|\partial_x^j\partial_\xi^k a(x,\xi)\bigr|\]

Given a symbol $a\in S_{\rho,\delta}^m(\mathbb R^{2d})$, its standard, or Kohn--Nirenberg, quantization is the operator initially defined on the Schwartz space $\mathcal S(\mathbb R^d)$ by 
\begin{equation*} \operatorname{Op}_\mathrm{kn}(a)\phi(x) \coloneqq \frac{1}{(2\pi)^d} \iint_{\mathbb R^d\times\mathbb R^d} e^{i\xi\cdot( x-y)} a(x,\xi)\phi(y)\,dy\,d\xi. \label{eq:KN-quantization} \end{equation*} 
Whenever the integral is not absolutely convergent, it is understood as an oscillatory integral.
Another quantization associated with the symbol $a$ is the Weyl quantization 
\begin{equation*} 
\operatorname{Op}_{\mathrm{w}}(a)\phi(x) \coloneqq \frac{1}{(2\pi)^d} \iint_{\mathbb R^d\times\mathbb R^d} e^{i\xi\cdot (x-y)} a\left(\frac{x+y}{2},\xi\right) \phi(y)\,dy\,d\xi. \label{eq:Weyl-quantization} 
\end{equation*}

\subsection{Magnetic pseudodifferential operators}



Consider a regular magnetic field $B$ described by smooth closed 2-forms on $\R^d$ that have components of class $BC^\infty(\R^d)$ \cite[Definition 2.8]{IftimiePurice2011}, i.e.:
\begin{equation}\label{eq:magdef}
\begin{aligned}
    &B \coloneqq\frac{1}{2} \sum_{1\le k,h \le d} B_{kh} (x) dx_k \land   dx_h, \quad \mbox{and}  \\
    &\sup_{x\in \R^d} |\partial_x^j B_{kh}(x)| < \infty\quad\forall j\in\N^d_0,
\end{aligned} 
\end{equation}
where $B_{kh}=-B_{hk}$. For every $M\geq0$ let
\begin{equation*}
    \Vert B\Vert_M\coloneqq \sum_{0\leq |j|\leq M}\sup_{1\leq k,h\leq d}\sup_{x\in\R^d}|\partial_x^jB_{kh}(x)|.
\end{equation*}
For a regular magnetic field we work with the associated magnetic potential $A$ in the transversal gauge:
\[A_h(x)= \sum_{1\le k \le d} \int_0^1 sx_k B_{kh } (sx) ds \]
Note the components are of class $C_{\mathrm{pol}}^\infty(\R^d)$, corresponding to functions in $C^\infty$ with polynomially bounded derivatives. Define also the uni-modular functions $\Lambda^A\colon \R^{2d} \to \C, \Omega^B\colon \R^{3d} \to \C$ with orientation given by the order of the points:
\begin{equation*}
\begin{aligned}
    \Lambda^A(x,y)&=\exp{[-i\cdot (\mbox{circulation of } A \mbox{ along the line with endpoints } x,y )]}\\
    &=\exp{[-i\cdot (\mbox{flux of } B \mbox{ trough the triangle with vertices } 0,x,y )]},
\end{aligned}
\end{equation*}
\begin{equation*}
\Omega^B(x,y,z)=\exp{[-i\cdot (\mbox{flux of } B \mbox{ trough the triangle with vertices } x,y,z )]}.
\end{equation*}
Then the magnetic pseudodifferential operator associated with the symbol $a$ is
\begin{equation*} 
\operatorname{Op}_\mathrm{w}^{A}(a)\phi(x) \coloneqq \frac{1}{(2\pi)^d} \iint_{\mathbb R^d\times\mathbb R^d} \Lambda^A(x,y)e^{i\xi\cdot (x-y)} a\left(\frac{x+y}{2},\xi\right) \phi(y)\,dy\,d\xi. \label{eq:Magnetic_PDO} 
\end{equation*}
In the absence of a magnetic field ($B\equiv0$), we have $\Lambda^0\equiv \Omega^0 \equiv 1$, therefore $\operatorname{Op}_\mathrm{w}^0(a)\equiv \operatorname{Op}_\mathrm{w}(a)$ for any symbol $a$.

\section{Auxiliary Results}\label{sec:aux}
We collect three ingredients: the Stokes factorization of $\Lambda^A$, derivative bounds for the flux factor, and a vector-valued Schur test.

Firstly, an application of Stokes Theorem gives Equation (3.3) in \cite{Cornean2024}: 
\begin{equation}\label{eq:Lambda_decomposition}
    \Lambda^A(x,y) = \Lambda^A(x,\alpha)\Lambda^A(\beta,y)\Lambda^A(\alpha,\beta) \Omega^B(\alpha ,x ,y) \Omega^B(\alpha , y, \beta )
\end{equation}
Secondly, a straightforward calculation gives Inequality (3.4) in \cite{Cornean2024}:
\begin{equation}\label{eq:Omega_boundedness}
\begin{aligned}
    |\partial_x^j \partial_y^k \partial_z^h &\Omega ^B ( x,y,z) | \\
    &\le C (1+\Vert B\Vert_{M})^{|j+k+h|} [1+ \mbox{diameter}(\triangle_{x,y,z} ) + \mbox{area}(\triangle_{x,y,z}) ]^{|j+k+h|}
\end{aligned}
\end{equation}
where $C>0$ and $M\in\N$ depends on $j,k,h$. Note we have specified the dependence on $B$ explicitly.




Lastly, we give a Schur test. 

\begin{proposition}\label{pro:gen_schur_test}
    Let $X$ be a Banach space, let $1 \leq p \leq \infty$ and let $ (A_{\alpha,\beta})_{\alpha,\beta\in\mathbb Z^d} \subseteq \mathcal L(X) $ be a family of bounded linear operators on $X$. 
    Assume that \[ C_1 \coloneqq \sup_{\alpha\in\mathbb Z^d} \sum_{\beta\in\mathbb Z^d} \|A_{\alpha,\beta}\|_{\mathcal L(X)} <\infty,\quad C_2 \coloneqq \sup_{\beta\in\mathbb Z^d} \sum_{\alpha\in\mathbb Z^d} \|A_{\alpha,\beta}\|_{\mathcal L(X)} <\infty.\]
    Then there exists a bounded linear operator $T\colon \ell^p(\mathbb Z^d;X) \longrightarrow \ell^p(\mathbb Z^d;X)$ satisfying
    \begin{equation}\label{eq:schur_op}
        (Tu)_\alpha \coloneqq \sum_{\beta\in\mathbb Z^d} A_{\alpha,\beta}u_\beta 
    \end{equation}
    for $u=(u_\beta)_{\beta\in\mathbb Z^d}\in\ell^p(\mathbb Z^d;X)$, and \[ \|T\|_{\mathcal L(\ell^p(\mathbb Z^d;X))} \leq C_1^{1/p'}C_2^{1/p}, \] where $p'$ is the conjugate exponent of $p$ with the convention $C_i^{1/\infty}=1$.
\end{proposition}

\begin{proof}
    Set $a_{\alpha,\beta}\coloneqq\Vert A_{\alpha,\beta}\Vert_{\mathcal{L}(X)}$. Note that by hypothesis $a_{\alpha,\cdot}=(a_{\alpha,\beta})_\beta\in\ell^1(\Z^d;\R)$ with uniform norm in $\alpha$.
    
    Let $u\in\ell^p(\Z^d;X)$. Using H\"older's inequality we get
    \begin{align*}
        \sum_{\beta}\Vert A_{\alpha,\beta}u_\beta\Vert_X\leq\sum_{\beta}a_{\alpha,\beta}\Vert u_\beta\Vert_X\leq\Vert a_{\alpha,\cdot}\Vert_{\ell^{p'}(\Z^d;\R)} \Vert u\Vert_{\ell^p(\Z^d;X)}<\infty,
    \end{align*}
    so the series $\sum_{\beta\in\mathbb Z^d} A_{\alpha,\beta}u_\beta$ converges in $X$ for every $\alpha$ and with a uniformly bounded norm in $\alpha$. This implies that Equation \eqref{eq:schur_op} gives a well-defined linear operator $T$ from $\ell^p(\Z^d;X)$ to $\ell^\infty(\Z^d;X)$. 
    
    Now for boundedness of $T$: if $1<p<\infty$, then
    \begin{align*}
        \Bigg(\sum_{\alpha}\Vert (Tu)_\alpha\Vert_X^p\Bigg)^{1/p}&\leq\Bigg(\sum_{\alpha}\Bigg( \sum_{\beta}a_{\alpha,\beta}^{\frac{1}{p'}}a_{\alpha,\beta}^\frac{1}{p}\Vert u_\beta\Vert_X\Bigg)^p\Bigg)^{1/p}\\
        &\leq\Bigg(\sum_{\alpha}\Bigg( \sum_{\beta}a_{\alpha,\beta}\Bigg)^{p-1} \sum_{\beta}a_{\alpha,\beta}\Vert u_\beta\Vert_X^p\Bigg)^{1/p}\\
        &\leq C_1^{1/p'}C_2^{1/p}\Vert u\Vert_{\ell^p(\Z^d;X)}.
    \end{align*}
    If $p=1$, then
    \begin{align*}
        \sum_{\alpha}\Vert (Tu)_\alpha\Vert_X&\leq\sum_{\alpha} \sum_{\beta}a_{\alpha,\beta}\Vert u_\beta\Vert_X\leq C_2\Vert u\Vert_{\ell^p(\Z^d;X)}.
    \end{align*}
    If $p=\infty$, then
    \begin{equation*}
        \sup_{\alpha}\Vert (Tu)_\alpha\Vert_X\leq\sup_{\alpha} \sum_{\beta}a_{\alpha,\beta}\Vert u_\beta\Vert_X\leq C_1\Vert u\Vert_{\ell^p(\Z^d;X)}.
    \end{equation*}
    This finishes the proof.
\end{proof}

\section{Proof of the main theorem}\label{sec:proof}






We first give two lemmata, one on the off-diagonal smoothness and decay of the \say{magnetic} twisted integral kernel and the other on the operator $\operatorname{Op}_\mathrm{w}^A (a)$ restricted to cubes. The first lemma is based on Proposition 2.1 in \cite{Taylor2023}.

\begin{lemma}\label{lem:off}
    Let $a\in S_{\rho,\delta}^m(\R^{2d})$ with $0<\rho\leq1$, $0\leq\delta\leq1$ and $m\in\R$. 
    
    The kernel $K^A(a)$ of $\mathrm{Op}_\mathrm{w}^A(a)$ is smooth off-diagonally and for every $R,N>0$ there exists $C_{R,N}>0$ and $N_{N,m,\rho,d}'\in\N$ independent of $a$ and $B$ such that
    \[|K^A(a)(x,y)|\leq\frac{C_{R,N}\Vert a\Vert_{S_{\rho,\delta}^m,N_{N,m,\rho,d}'}}{|x-y|^N}\]
    when $|x-y|\geq R$.
\end{lemma}

\begin{proof}
    For any $N>0$ large enough so that the function $\langle\xi\rangle^{m-2\rho N}$ is absolutely integrable we have that 
    \[
    \left| \triangle_\xi^N a\left(\frac{x+y}{2}, \xi\right) \right| \le C_N \langle \xi \rangle ^{m-2\rho N}.
    \]
    So, the following integral is well-defined and is uniformly bounded:
    \[
    \left|\int_{\R^d} e^{i\xi\cdot (x-y)}  (-\triangle_\xi)^N\left[a \left( \frac{x+y}{2} ,\xi \right) \right] d\xi \right| \le C_N\Vert a\Vert_{S_{\rho,\delta}^m,2N}
    \]
    It follows that the function $f_N(x,y)\coloneqq\int_{\R^d} e^{i\xi\cdot (x-y)}  (-\triangle_\xi)^N\left[a \left( \frac{x+y}{2} ,\xi \right) \right] d\xi$ is differentiable of some order $k_{N,\delta}$ that goes to infinity with $N$. Additionally, we see that
    \[
    \begin{split}
     \frac{1}{(2\pi)^d}\Lambda^A(x,y)f_N(x,y)&= |x-y|^{2N} K^A(a)(x,y)
    \end{split}
    \]
    as distributions when $|x-y|>0$, which shows that $K^A(a)$ is smooth off-diagonally and gives the estimate.
\end{proof}

\begin{lemma}\label{lem:cube}
    Let $a\in S_{\rho,\delta}^m(\R^{2d})$ with $0<\rho\leq1$, $0\leq\delta\leq1$ and $m\in\R$. Let $K^A(a)$ be the kernel of $\mathrm{Op}_\mathrm{w}^A(a)$. Lastly fix $p\in[1,\infty)$.

    For two functions $g_\alpha,g_\beta\in C^\infty(\R^d)$ with $\alpha,\beta\in\Z^d$, $\supp g_\alpha\subseteq[-1,1]^d+\alpha$, $\supp g_\beta\subseteq[-1,1]^d+\beta$, and $0\leq g_\alpha,g_\beta\leq1$, the operator $T$ with distributional kernel $g_\alpha(x)K^A(x,y)g_\beta(y)$ is bounded on $L^p(\R^d)$ if $\mathrm{Op}_\mathrm{w}(a)$ is bounded. Moreover, for every $R>0$, there exists $C_R>0$ and $N,M\in\N$ independent of $a$ and $B$ such that
    \[\Vert T\Vert_{\mathcal{L}(L^p)}\leq C_R(1+\Vert B\Vert_M)^N\Vert\mathrm{Op}_\mathrm{w}(a)\Vert_{\mathcal{L}(L^p)}\]
    with $r=|\alpha-\beta|\leq R$.
\end{lemma}

\begin{proof}
    Equation \eqref{eq:Lambda_decomposition} allows one to write 
    \begin{align*}
        &g_\alpha(x)K^A(x,y)g_\beta(y) \\
        &=g_\alpha(x) g_\beta (y) \Lambda^A(x,\alpha)\Lambda^A(\beta,y)\Lambda^A (\alpha,\beta) \Omega^B(\alpha,x,y) \Omega^B(\alpha,y,\beta) K^0(a)(x,y).
    \end{align*}
    Construct a function $h_\alpha\in C^\infty(\R^d)$ such that $h_\alpha g_\alpha=g_\alpha$, $\supp h_\alpha\subseteq [-2,2]^d+\alpha$, and $0\leq h_\alpha\leq1$, and construct $h_\beta$ with similar properties but for $g_\beta$.
    Define
    \[
    t(x,y) \coloneqq  h_\alpha(x+\alpha)h_\beta(y+\beta) \Lambda^A (\alpha,\beta) \Omega^B(\alpha,x+\alpha,y+\beta) \Omega^B(\alpha,y+\beta,\beta),
    \]
    so that
    \begin{equation}\label{eq:K_rewriting}
        g_\alpha(x)K^A(x,y)g_\beta(y) = g_\alpha(x) g_\beta(y) \Lambda^A(x,\alpha)\Lambda^A(\beta,y)t(x-\alpha,y-\beta) K^0(a)(x,y).
    \end{equation}
    The function $t$ has support in $[-2,2]^{2d}\subset [-\pi,\pi]^{2d}$ and considering $t$ as an element of $C^\infty([-\pi,\pi]^{2d})$, we can expand it in a Fourier series: 
    \begin{equation*}
        t(x,y) = \sum_{n,m\in\Z^d} c_{n,m} e^{in\cdot x} e^{im\cdot y}
    \end{equation*}
    with 
    \[
    c_{n, m} =\frac{1}{(2\pi)^{2d}} \int_{[-\pi,\pi]^{2d}} dx dy  e^{-in\cdot x} e^{-im\cdot y} t(x,y).
    \]
    Due to the smoothness of $t$ and the fact that $t$ vanishes at the boundary we conclude that the coefficients $c_{n, m}$ decay super-polynomially and the Fourier series converges absolutely. Furthermore, the supports of $h_\alpha,h_\beta$ and size of $r=|\alpha-\beta|$ controls the diameters and areas of the triangles with vertices $\alpha,x+\alpha,y+\beta$ and $\alpha,y+\beta,\beta$, which implies by Equation \eqref{eq:Omega_boundedness} that there exists a constant $C_R$, depending only on $R$, and an $r$-independent constant $M\in\N$ such that for $r\leq R$, 
    \begin{equation*}
        |\partial_x^j\partial_y^k t (x,y)|\le C_R^{|j+k|}(1+\Vert B\Vert_{M})^{|j+k|}.
    \end{equation*}
    Now this implies that, for $N>d/2$,
    \[
    |c_{n, m} |\le \frac{C_R^{4N}(1+\Vert B\Vert_{M'})^{4N}}{[(1+|n|^2)(1+|m|^2)]^N}.
    \]
    Returning to Equation \eqref{eq:K_rewriting}, the presence of the localizing functions $g_\alpha$ and $g_\beta$ ensures the following identity holds in the distributional sense:
    \begin{align*}
        g_\alpha(x)&K^A(x,y)g_\beta(y) \\
        &=\sum_{n,m\in\Z^d} c_{n,m} g_\alpha(x) e^{in\cdot (x-\alpha)}\Lambda^A(x,\alpha)g_\beta (y) e^{im\cdot (y-\beta)} \Lambda^A(\beta,y) K^0(a)(x,y).
    \end{align*}
   
    Under the hypothesis that $\operatorname{Op}_\mathrm{w}(a) \colon\mathcal{S}(\R^d) \to \mathcal{S}(\R^d)$ is $\Lp{p}$-bounded, we have that for any $\phi \in \mathcal{S}(\R^d)$, $\operatorname{Op}_\mathrm{w}(a)\phi$ is an $\Lp{p}$-function and the same holds for $T\phi$ by the above. In particular the following estimate holds as the function $\Lambda^A$ is uni-modular:
    \[
    \begin{split}
        \|T&\phi\|_{\Lp{p}}   \le \sum_{n,m\in\Z^d} |c_{n,m} | \\ &\quad \times \left\|g_\alpha(x) e^{in\cdot (x-\alpha)}\Lambda^A(x,\alpha)\operatorname{Op}_\mathrm{w}(a)\left[g_\beta (y) e^{im\cdot (y-\beta)} \Lambda^A(\beta,y)\phi(y)\right] (x) \right\|_{\Lp{p}} \\
        &\le \sum_{n,m\in\Z^d} |c_{n,m} | \left\|\operatorname{Op}_\mathrm{w}(a)\left[g_\beta (y) e^{im\cdot (y-\beta)} \Lambda^A(\beta,y)\phi(y)\right] (x) \right\|_{\Lp{p}} \\
        &\le \sum_{n,m\in\Z^d} |c_{n,m} | \left\|\operatorname{Op}_\mathrm{w}(a) \right\|_{\mathcal{L}(\Lp{p})} \left\|g_\beta (y) e^{im\cdot (y-\beta)} \Lambda^A(\beta,y)\phi(y)  \right\|_{\Lp{p}} \\
        & \le \sum_{n,m\in\Z^d} |c_{n,m} | \left\|\operatorname{Op}_\mathrm{w}(a) \right\|_{\mathcal{L}(\Lp{p})} \left\|\phi  \right\|_{\Lp{p}} \\
        &\le C_R'(1+\Vert B\Vert_{M'})^{N'} \left\|\operatorname{Op}_\mathrm{w}(a) \right\|_{\mathcal{L}(\Lp{p})} \left\|\phi  \right\|_{\Lp{p}}
    \end{split}
    \]
    for some $C_R'>0$ and $N'\in\N$. Since $\mathcal{S}(\R^d)$ is dense in $\Lp{p}(\R^d)$ we have
    \begin{equation*}\label{eq:norm_estimate_Op_cellular}
        \|T\|_{\mathcal{L}(\Lp{p})} \le C_R'(1+\Vert B\Vert_{M'})^{N'} \|\operatorname{Op}_\mathrm{w}(a)\|_{\mathcal{L}(\Lp{p})}.\qedhere
    \end{equation*}
\end{proof}

Now we move on to the proof of Theorem \ref{thm:mainresult}. The proof for $p=\infty$ follows from the case $p=1$ by duality arguments, so we first prove the case $p<\infty$ and leave $p=\infty$ for the very end.

\subsection{When $p<\infty$}

\begin{proof}[Proof of Theorem \ref{thm:mainresult} when $p\in[1,\infty)$]
    Fix $p\in[1,\infty)$. 
    
    Let $K^A(a),K^0(a)\in \mathcal{S}'(\R^{2d})$ be the kernels of $\mathrm{Op}_\mathrm{w}^A(a),\mathrm{Op}_\mathrm{w}(a)$, and note that $K^A(a)(x,y)=\Lambda^A(x,y)K^0(a)(x,y)$ as distributions. The two implications in Theorem \ref{thm:mainresult} can be proven using the exact same methods since \[K^0(a)(x,y)=\overline{\Lambda^A(x,y)}K^A(a)(x,y)=\Lambda^A(y,x)K^A(a)(x,y),\] so we only prove that if $\operatorname{Op}_\mathrm{w} (a)$ is bounded, then $\operatorname{Op}_\mathrm{w}^A(a)$ is bounded.
    
    Construct a non-negative function $g \in C^\infty (\R^d)$ such that \begin{equation*}\label{eq:g}
        \supp g \subset [-1,1]^d, \quad \mbox{ and } \sum_{\alpha\in\Z^d} g(x-\alpha)=1 \quad\mbox{pointwise.}
    \end{equation*}
    Now for $\alpha,\beta \in \Z^d$, set $g_\alpha(x) = g(x-\alpha)$ and set \[K^A(a)_{\alpha,\beta}(x,y) = g_\alpha(x)K^A(a)(x,y)g_\beta(y)\] with the associated operator being $\operatorname{Op}_\mathrm{w}^A(a)_{\alpha,\beta}.$
    The properties of $g$ ensures that for any $\phi\in \mathcal{S}(\R^d)$ we have that $\sum_{\alpha\in\Z^d}g_\alpha \phi$ converges to $\phi$ absolutely in $\mathcal{S}(\R^d)$: for $\gamma\in\N_0^d$ and $N\in\N$ we have
    \begin{align*}
        \sum_{\alpha\in\Z^d}&\sup_{x\in\R^d}|\langle x\rangle^N\partial^\gamma(g_\alpha(x)\phi(x))|\\
        &\leq2^{\frac{d+1}{2}}\sum_{\alpha\in\Z^d}\langle \alpha\rangle^{-d-1}\sup_{x\in\R^d}|\langle x\rangle^{N+d+1}\langle x-\alpha\rangle^{d+1}\partial^\gamma(g_\alpha(x)\phi(x))|\\
        &\leq C\sum_{\gamma'\leq\gamma}\sup_{x\in\R^d}|\langle x\rangle^{N+d+1}\partial^{\gamma'}\phi(x)|\sum_{\alpha\in\Z^d}\langle \alpha\rangle^{-d-1}<\infty
    \end{align*}
    and $\sum_{\alpha\in\Z^d}g_\alpha \phi=\phi$ holds pointwise.
    So for any test function $\phi\in \mathcal{S}(\R^d)$, the following series converges in $\mathcal{S}(\R^d)$: $$
    \sum_{\alpha,\beta \in \Z^d} \operatorname{Op}_\mathrm{w}^A(a)_{\alpha,\beta} \phi=\operatorname{Op}_\mathrm{w}^A(a)\phi $$
    Now fix $R>4d$ and define the diagonal and off-diagonal operators $$\operatorname{Op}_\mathrm{w}^A(a)_\mathrm{D} = \sum_{|\alpha-\beta| <R} \operatorname{Op}_\mathrm{w}^A(a)_{\alpha,\beta}, \quad\operatorname{Op}_\mathrm{w}^A(a)_\mathrm{OD} = \sum_{|\alpha-\beta| \ge R } \operatorname{Op}_\mathrm{w}^A(a)_{\alpha,\beta}.  $$
    For any $\phi \in \mathcal{S}(\R^d)$ we have $\operatorname{Op}_\mathrm{w}^A(a) \phi = \operatorname{Op}_\mathrm{w}^A(a)_\mathrm{D} \phi + \operatorname{Op}_\mathrm{w}^A(a)_\mathrm{OD} \phi$. The off-diagonal term $\operatorname{Op}_\mathrm{w}^A(a)_\mathrm{OD}$ is $L^p$-bounded by Lemma \ref{lem:off} and the scalar Schur test \cite{Grafakos}, therefore if we can prove the $L^p$-boundedness of the diagonal part $\operatorname{Op}_\mathrm{w}^A(a)_\mathrm{D}$, then the proof is completed.
    
    So let us focus on the diagonal part $\operatorname{Op}_\mathrm{w}^A(a)_\mathrm{D}$. It follows from Lemma \ref{lem:cube} and the hypothesis that each $\operatorname{Op}_\mathrm{w}^A(a)_{\alpha,\beta}$ is $L^p$-bounded and the operator norms are uniformly bounded for $|\alpha-\beta|<R$. To finish the proof we use the vector-valued Schur test (Proposition \ref{pro:gen_schur_test}).
    
    Consider the space of functions $\mathcal{B}^p=\ell^p (\Z^d;\Lp{p}([0,1)^d))$.
    It is a Banach space with the norm 
    \[
    \|F\|_{\mathcal{B}^p} = \left( \sum_{\alpha\in\Z^d} \|F(\alpha)\|_{\Lp{p}([0,1)^d)}^p\right)^{\frac{1}{p}}.
    \]
    Define the linear operators $\mu\colon\Lp{p}(\R^d) \to \mathcal{B}^p$, $\nu\colon \mathcal{B}^p \to \Lp{p}(\R^d)$ as follows:
    \[\mu(f)(\alpha) |_x=f(x+\alpha)\quad \forall\alpha\in\Z^d, x\in [0,1)^d,\]
    \[\nu(F)(x) = F(\alpha)|_{x-\alpha}  \quad \mbox{for the unique } \alpha \in \Z^d \ | \ x\ \in [0,1)^d+\alpha .\]
    The function $\mu$ is an isometry since it holds that
    \[
    \begin{split}
        \|f\|_{\Lp{p}} &= \left( \int_{\R^d} dx |f(x)|^p\right)^{\frac{1}{p}} = \left( \sum_{\alpha\in\Z^d} \int_{[0,1)^d-\alpha} dx |f(x)|^p\right)^{\frac{1}{p}} \\ 
        & = \left\{ \sum_{\alpha\in\Z^d}\left[ \left( \int_{[0,1)^d}  dx |f(x+\alpha)|^p\right)^\frac{1}{p}\right]^p\right\}^{\frac{1}{p}} = \|\mu(f)\|_{\mathcal{B}^p}.
    \end{split}
    \]
    An identical proof shows that $\nu$ is also an isometry with $\mu\circ \nu = \mathrm{Id}_{\mathcal{B}^p}, \nu\circ \mu = \mathrm{Id}_{\Lp{p}}$.
    Define the families of operators $\{\mu_\alpha\}_{\{\alpha\in\Z^d\}} \colon L^p(\R^d)\to L^p([0,1)^d) $ and $\{\nu_\beta\}_{\{\beta\in\Z\}} \colon L^p([0,1)^d)\to L^p(\R^d)$ such that for all $\alpha,\beta\in \Z^d, x \in \R^d, y\in [0,1)^d$ the operators act as follows:
    \[
    \mu_\alpha(f)(x)=f(x+\alpha)\quad \mbox{and} \quad \nu_\beta (f) (y) = \begin{cases}
        f(y-\beta) &\mbox{if } y-\beta \in[0,1)^d, \\0&\mbox{otherwise.}
    \end{cases}
    \]
    Then the operators $$\widetilde{\operatorname{Op}_\mathrm{w}^A}(a)_{\alpha,\beta,\alpha',\beta'} = \mu_{\alpha'} \circ \operatorname{Op}_\mathrm{w}^A(a)_{\alpha,\beta} \circ \nu_{\beta'} \colon\Lp{p}([0,1)^d)\to \Lp{p}([0,1)^d)$$
    are bounded with $$\left\|\widetilde{\operatorname{Op}_\mathrm{w}^A}(a)_{\alpha,\beta,\alpha',\beta'}\right\|_{\mathcal{L}(\Lp{p}([0,1)^d))}\le C(1+\Vert B\Vert_{M})^{N}\|\operatorname{Op}_\mathrm{w} (a)\|_{\mathcal{L}(L^p)}$$
    for some $C>0$ and $M,N\in\N$ independent of $a$ and $B$ from Lemma \ref{lem:cube} and since all $\mu_{\alpha'}$ and $\nu_{\beta'}$ have operator norms equal to one.
    Most importantly, they decompose the action of the operator $\operatorname{Op}_\mathrm{w}^A(a)_{\alpha,\beta}$ as follows:
    \begin{equation}\label{eq:decomp}
        \left[\mu\circ \operatorname{Op}_\mathrm{w}^A(a)_{\alpha,\beta} \circ \nu \right](F)(\alpha') = \sum_{\beta'\in\Z^d} \widetilde{\operatorname{Op}_\mathrm{w}^A}(a)_{\alpha,\beta,\alpha',\beta'}F(\beta').
    \end{equation}

    Now, the supports of the functions $g_\alpha,g_\beta$ imply that
    \[
    \supp \left( \operatorname{Op}_\mathrm{w}^A(a)_{\alpha,\beta} (\phi )\right)\subset [-1,1]^d+\alpha \quad \forall\phi\in \mathcal{S }(\R^d), \quad \alpha, \beta\in\Z^d, \mbox{ and}
    \]
    \[
    \operatorname{Op}_\mathrm{w}^A(a)_{\alpha,\beta} (\phi )\equiv 0 \quad \mbox{ if} \quad \supp (\phi) \cap ( [-1,1]^d +\beta) = \emptyset \quad \forall\alpha,\beta \in\Z^d.
    \]
    So in the series in Equation \eqref{eq:decomp} there is only a finite number of $\widetilde{\operatorname{Op}_\mathrm{w}^A}(a)_{\alpha,\beta,\alpha',\beta'}$ different than zero at a time and they are all concentrated around $(\alpha,\beta)\in\Z^{2d}$: a rough estimate would be \begin{equation}\label{eq:when_blocks_are_zero}
    \widetilde{\operatorname{Op}_\mathrm{w}^A}(a)_{\alpha,\beta,\alpha',\beta'} = 0\quad  \mbox{if}\quad \sum_{n=1}^d |\alpha_n-\alpha_n'|+ |\beta_n-\beta_n'|>2d.
    \end{equation}
    For this reason, $$\sum_{|\alpha-\beta|<R} \widetilde{\operatorname{Op}_\mathrm{w}^A}(a)_{\alpha,\beta,\alpha',\beta'} \eqqcolon\widetilde{\operatorname{Op}_\mathrm{w}^A}(a)_{\mathrm{D},\alpha',\beta'}$$
    is also well-defined for every $\alpha',\beta'\in\Z^d$ as there is a finite number of non-zero entries on each side.
    Using also Equation \eqref{eq:when_blocks_are_zero} we conclude that 
    \begin{align*}\Big\|\widetilde{\operatorname{Op}_\mathrm{w}^A}&(a)_{\mathrm{D},\alpha',\beta'}\Big\|_{\mathcal{L}(L^p([0,1)^d))}\\
    &\le \begin{cases}
        C'(1+\Vert B\Vert_{M})^{N}\|\operatorname{Op}_\mathrm{w}(a)\|_{\mathcal{L}(L^p(\R^d))} &\mbox{if } |\alpha'-\beta'|<R+2d, \\ 0 & \mbox{otherwise},
    \end{cases}\end{align*}
    for $C'>0$ and $M,N\in\N$.
    This allows for the use of the vector-valued Schur test (Proposition \ref{pro:gen_schur_test}), implying that the operator
    \[
    \widetilde{\operatorname{Op}_\mathrm{w}^A}(a)_{\mathrm{D}}(F)(\alpha')\coloneqq\sum_{\beta'\in\Z^d}\widetilde{\operatorname{Op}_\mathrm{w}^A}(a)_{\mathrm{D},\alpha',\beta'} F(\beta')
    \]
    is bounded. However, it is also true that $$\widetilde{\operatorname{Op}_\mathrm{w}^A}(a)_\mathrm{D} = \mu\circ \operatorname{Op}_\mathrm{w}^A (a)_\mathrm{D}\circ \nu,$$ and the boundedness of the right hand side follows.
\end{proof}

\begin{remark}\label{rem:bound}
    Note that since the estimate on the off-diagonal part is $B$-independent (only $|\Lambda^A|=1$ is used), the above proof shows that
    \begin{align*}
        \lVert\mathrm{Op}_\mathrm{w}^A(a)_\mathrm{D}\rVert_{\mathcal{L}(L^p)}&\leq C(1+\Vert B\Vert_M)^{K}\Vert\mathrm{Op}_\mathrm{w}(a)\Vert_{\mathcal{L}(L^p)}\\
        \lVert\mathrm{Op}_\mathrm{w}^A(a)_\mathrm{OD}\rVert_{\mathcal{L}(L^p)}&\leq C'\Vert a\Vert_{S_{\rho,\delta}^m,N}
    \end{align*}
    for some $C,C'>0$ and $K,N,M\in\N$ independent of $a,B$. So
    \[\lVert\mathrm{Op}_\mathrm{w}^A(a)\rVert_{\mathcal{L}(L^p)}\leq C''(1+\Vert B\Vert_M)^{K'}(\Vert a\Vert_{S_{\rho,\delta}^m,N'}+\Vert\mathrm{Op}_\mathrm{w}(a)\Vert_{\mathcal{L}(L^p)})\]
    for $C''>0$ and $N',M,K'\in\N$ independent of $a,B$.
\end{remark}

\subsection{When $p=\infty$}

\begin{proof}[Proof of Theorem \ref{thm:mainresult} when $p=\infty$]
    Just as in the proof of the case $p<\infty$, we only prove that if $\mathrm{Op}_\mathrm{w}(a)$ is $L^\infty$-bounded, then $\mathrm{Op}_\mathrm{w}^A(a)$ is $L^\infty$-bounded.

    From Lemma 1.2 in \cite{Purice2007} and Proposition 21 in \cite{Purice2004} we know that $\operatorname{Op}_\mathrm{w}^A(a)$ is extendable to a continuous operator on $\mathcal{S}'(\R^d)$, so it is defined on $L^\infty(\R^d)$. Considering the duality bracket $\langle\cdot, \cdot\rangle$ for $\mathcal{S}'(\R^d)\times \mathcal{S}(\R^d)$, we get the identity 
    \[
    \left\langle \operatorname{Op}_\mathrm{w}^A(a) f, \phi \right\rangle = \left\langle f , \overline{\operatorname{Op}_\mathrm{w}^A(\overline{a}) \overline{\phi}} \right\rangle
    \]
    for $f\in L^\infty(\R^d)$ and $\phi\in\mathcal{S}(\R^d)$.

    Suppose $\mathrm{Op}_\mathrm{w}(a)$ is $L^\infty$-bounded. Then
    \begin{equation*}
        \begin{aligned}
            \Vert\operatorname{Op}_\mathrm{w}(a)\Vert_{\mathcal{L}(L^\infty)}&=\sup_{\Vert f\Vert_{L^\infty}\leq1,\Vert\phi\Vert_{L^1}\leq1,\phi\in\mathcal{S}(\R^d)}\left|\left\langle \operatorname{Op}_\mathrm{w}(a) f, \phi \right\rangle\right|\\
            &=\sup_{\Vert\phi\Vert_{L^1}\leq1,\phi\in\mathcal{S}(\R^d)}\Vert\operatorname{Op}_\mathrm{w}(\overline{a})\overline{\phi} \Vert_{L^1}
        \end{aligned}
    \end{equation*}
    and so $\mathrm{Op}_\mathrm{w}(\overline{a})$ is $L^1$-bounded. This implies that $\mathrm{Op}_\mathrm{w}^A(\overline{a})$ is $L^1$-bounded (Theorem \ref{thm:mainresult} $p=1$) and for a given $f\in L^\infty(\R^d)$ we have
    \begin{equation*}
        \begin{aligned}
            \sup_{\Vert\phi\Vert_{L^1}\leq1,\phi\in\mathcal{S}(\R^d)}\left|\left\langle \operatorname{Op}_\mathrm{w}^A(a)f,\phi \right\rangle\right|&=\sup_{\Vert\phi\Vert_{L^1}\leq1,\phi\in\mathcal{S}(\R^d)}\left|\left\langle f, \overline{\operatorname{Op}_\mathrm{w}^A(\overline{a})\overline{\phi}} \right\rangle\right|\\
            &\leq\Vert f\Vert_{L^\infty}\Vert\mathrm{Op}_\mathrm{w}^A(\overline{a})\Vert_{\mathcal{L}(L^1)},
        \end{aligned}
    \end{equation*}
    which shows that $\operatorname{Op}_\mathrm{w}^A(a)$ leaves $L^\infty(\R^d)$ invariant. Furthermore,
    \begin{equation*}
        \begin{aligned}
            \sup_{\Vert f\Vert_{L^\infty}\leq1}\Vert\operatorname{Op}_\mathrm{w}^A(a)f\Vert_{L^\infty}&=\sup_{\Vert f\Vert_{L^\infty}\leq1,\Vert\phi\Vert_{L^1}\leq1,\phi\in\mathcal{S}(\R^d)}\left|\left\langle f, \overline{\operatorname{Op}_\mathrm{w}^A(\overline{a})\overline{\phi}} \right\rangle\right|\\
            &=\Vert\mathrm{Op}_\mathrm{w}^A(\overline{a})\Vert_{\mathcal{L}(L^1)}<\infty,
        \end{aligned}
    \end{equation*}
    so $\operatorname{Op}_\mathrm{w}^A(a)$ is $L^\infty$-bounded.
\end{proof}

\section{Proof of Corollary \ref{cor:mainresult}}\label{sec:appl}

We now combine Theorem \ref{thm:mainresult} with known results. First, we have an adapted theorem on the change of quantization transformation, see Theorem 4.5 in \cite{Hrmander1979} and Theorem 18.1.7 and the discussion before Theorem 18.1.35 in \cite{Hrmander1985}.

\begin{theorem}\label{thm:changequant}
    Let $m\in\R$ and $0<\rho,0\leq\delta\leq\rho\leq1,\delta<1$. Then $e^{\frac{i}{2}\nabla_x\cdot\nabla_\xi}$ is a homeomorphism on $S_{\rho,\delta}^m(\R^{2d})$ and
    \begin{equation}\label{eq:changequant}
        \mathrm{Op_{kn}}(a)=\mathrm{Op_w}(e^{\frac{i}{2}\nabla_x\cdot\nabla_\xi}a)
    \end{equation}
    for all $a\in S_{\rho,\delta}^m(\R^{2d})$.
\end{theorem}

\begin{proof}
    The proof that $e^{\frac{i}{2}\nabla_x\cdot\nabla_\xi}$ and its inverse $e^{-\frac{i}{2}\nabla_x\cdot\nabla_\xi}$ leave $S_{\rho,\delta}^m(\R^{2d})$ invariant and that Equation \eqref{eq:changequant} hold can be adapted from \cite{Hrmander1985}.

    The continuity statement for $e^{\frac{i}{2}\nabla_x\cdot\nabla_\xi}$ and its inverse follow from the Closed Graph Theorem: let $(a_n)_{n\in\N}$ be a sequence in $S_{\rho,\delta}^m(\R^{2d})$. Suppose that $a,b\in S_{\rho,\delta}^m(\R^{2d})$ exist such that $a_n\stackrel{n\rightarrow\infty}{\longrightarrow} a$ and $e^{\frac{i}{2}\nabla_x\cdot\nabla_\xi}a_n\stackrel{n\rightarrow\infty}{\longrightarrow} b$ in $S_{\rho,\delta}^m(\R^{2d})$. Then both limits hold in $\mathcal{S}'(\R^{2d})$, and moreover, $e^{\frac{i}{2}\nabla_x\cdot\nabla_\xi}$ is continuous on $\mathcal{S}'(\R^{2d})$, so $e^{\frac{i}{2}\nabla_x\cdot\nabla_\xi}a_n\stackrel{n\rightarrow\infty}{\longrightarrow} e^{\frac{i}{2}\nabla_x\cdot\nabla_\xi}a$ in $\mathcal{S}'(\R^{2d})$. Thus $b=e^{\frac{i}{2}\nabla_x\cdot\nabla_\xi}a$ and the proof is finished.
\end{proof}

Second, we have:
\begin{theorem}\label{thm:pboundedkn}(\cite[Theorem 4.1]{Somersalo1988})
    Let $p\in(1,\infty)$ and $0\leq\rho,0\leq\delta\leq\rho\leq1,\delta<1$. If $m \le -d (1-\rho) \left|\frac{1}{p} -\frac{1}{2} \right|$ and $a\in S_{\rho,\delta}^m(\R^{2d})$, then $\mathrm{Op_{kn}}(a)$ is bounded on $L^p (\R^d)$. Moreover, there exists $C>0$ and $N\in\N$ independent of $a$ such that
    \[\Vert\mathrm{Op_{kn}}(a)\Vert_{\mathcal{L}(L^p)}\leq C\Vert a\Vert_{S_{\rho,\delta}^m,N}.\]
\end{theorem}
\begin{proof}
The only thing not stated explicitly in \cite{Somersalo1988} and \cite[pp. 322--323]{Stein2016} is the bound on the operator norm. But one can simply use the Closed Graph Theorem as done in the proof of Lemma 2.6 in \cite{Benyi2013}:

Let $(a_n)_{n\in\N}$ be a sequence in $S_{\rho,\delta}^m(\R^{2d})$, and suppose that $(a_n)_{n\in\N}$ converges to $a\in S_{\rho,\delta}^m(\R^{2d})$ and that $(\mathrm{Op}_\mathrm{kn}(a_n))_{n\in\N}$ converges to $T\in\mathcal{L}(L^p(\R^d))$. Then $a_n\stackrel{n\rightarrow\infty}{\longrightarrow} a$ in $\mathcal{S}'(\R^{2d})$ implying that \[\mathrm{Op}_\mathrm{kn}(a_n)\stackrel{n\rightarrow\infty}{\longrightarrow}\mathrm{Op}_\mathrm{kn}(a)\] in $\mathcal{L}(\mathcal{S}(\R^d),\mathcal{S}'(\R^d))$. But $\mathrm{Op}_\mathrm{kn}(a_n)\stackrel{n\rightarrow\infty}{\longrightarrow}T$ in $\mathcal{L}(L^p(\R^d))\hookrightarrow\mathcal{L}(\mathcal{S}(\R^d),\mathcal{S}'(\R^d))$, so $T=\mathrm{Op}_\mathrm{kn}(a)$.

This proves that $S_{\rho,\delta}^m(\R^{2d})\ni a\to\mathrm{Op}_\mathrm{kn}(a)\in\mathcal{L}(L^p(\R^d))$ is a continuous mapping, which is equivalent with the existence of the stated bound on the operator norm.
\end{proof}

This now leads to:

\begin{proof}[Proof of Corollary \ref{cor:mainresult}]
    Let $p,\rho,\delta,m$ satisfy the criteria in the Corollary. 
    
    Let $a\in S_{\rho,\delta}^m(\R^{2d})$. From Theorem \ref{thm:changequant} and Theorem \ref{thm:pboundedkn} we see that $\mathrm{Op}_\mathrm{w}(a)$ is bounded on $L^p(\R^d)$ and that $S_{\rho,\delta}^m(\R^{2d})\ni b\to\mathrm{Op}_\mathrm{w}(b)\in\mathcal{L}(L^p(\R^d))$ is a continuous mapping. Now using Theorem \ref{thm:mainresult} we conclude that $\mathrm{Op}_\mathrm{w}^A(a)$ is bounded on $L^p(\R^d)$, and using Remark \ref{rem:bound} we get
    \begin{align*}
        \lVert\mathrm{Op}_\mathrm{w}^A(a)\rVert_{\mathcal{L}(L^p)}&\leq C(1+\Vert B\Vert_M)^{K}(\Vert a\Vert_{S_{\rho,\delta}^m,N}+\Vert\mathrm{Op}_\mathrm{w}(a)\Vert_{\mathcal{L}(L^p)})\\
        &\leq C'(1+\Vert B\Vert_M)^{K}\Vert a\Vert_{S_{\rho,\delta}^m,N'}
    \end{align*}
    for some $C,C'>0$ and $N,M,K,N'\in\N$ independent of $a,B$.
\end{proof}

\begin{remark}
    Corollary \ref{cor:mainresult} does not say anything on $L^p$-boundedness with $p=1$ or $p=\infty$ whereas Theorem \ref{thm:mainresult} does include these cases. But of course any results on the $L^p$-boundedness with $p=1$ or $p=\infty$ for $\mathrm{Op_w}(a)$, with $a$ satisfying the hypothesis of Theorem \ref{thm:mainresult}, also hold for $\mathrm{Op}_\mathrm{w}^A(a)$.
\end{remark}

\section{Open questions and extensions}\label{sec:open}
Corollary \ref{cor:mainresult} covers the classical result almost completely \cite{Somersalo1988} with the exception of the limiting case $\rho=0$. When $\rho=0$, the distributional kernels lack pseudo-locality \cite{Stein2016}, making the control of the off-diagonal term in our proof of Theorem \ref{thm:mainresult} impossible. To establish this case, one could prove boundedness from $L^\infty$ to $BMO$ as in \cite{Fefferman1973} and on Hardy spaces $h^p,H^p$ as in \cite{Somersalo1988}, and then interpolate with the magnetic Calder\'on--Vaillancourt Theorem \cite{Purice2007}. This approach would also in and of itself be interesting, as improvements of our result for $p=1, p=\infty$ are expected to hold only between pairs of spaces as $L^\infty \to BMO, L^1\to h^1$. Maybe it is also possible to obtain equivalent results for $p< 1$ using $ h^p,H^p,L^p$, but we expect completely different techniques to be needed.


An extension of our work could be to work with magnetic Fourier integral operators. $L^2$-boundedness (and boundedness in magnetic $L^2$-Sobolev spaces) has been proven in the magnetic case \cite{IftimiePurice2011}. We think this could be solved by proving an equivalence theorem as done in this paper.

It is also possible that one could apply our result to the magnetic pseudodifferential super calculus \cite{LeeLein2022}, e.g. proving boundedness of magnetic pseudodifferential super operators on the spaces $\mathcal{L}(L^{p}(\R^d),L^{q}(\R^d))$. For magnetic pseudodifferential super operators a Hilbert--Schmidt class boundedness result was proven in \cite{LeeLein2025} and boundedness for some Banach spaces of operators in \cite{CorneanThorn2026}.


\section*{Acknowledgements} 
The authors acknowledge support from the Independent Research Fund Denmark--Natural Sciences, grant DFF–5281-00046B, and from DNRF Centre CLASSIQUE sponsored by the Danish National Research Foundation, grant nr. 187.

The authors acknowledge fruitful discussions with Horia Cornean.

\section*{Declarations}
\subsection*{AI disclosure} 

Proposition \ref{pro:gen_schur_test} was formulated with the help of AI. 
\subsection*{Competing interests} 

The authors have no competing interests.

\bibliographystyle{amsalpha} 
\bibliography{references} 
\end{document}